\documentclass{amsart}
\usepackage{verbatim}

\newtheorem{theorem}{Theorem}
\newtheorem{lemma}[theorem]{Lemma}
\theoremstyle{definition}

\theoremstyle{remark}
\newtheorem{remark}[]{Remark}
\numberwithin{equation}{section} \theoremstyle{corollary}

\theoremstyle{proposition}

\newfont{\EUL}{eufm10 scaled 1000}

\newcommand\R{\mathbb{R}}

\newcommand\Ker{{\rm Ker}}

\newcommand\Lt{\mbox{\EUL t}}

\renewcommand\k{\mbox{\EUL k}}

\newcommand\m{\mbox{\EUL m}}

\newcommand\Ad{{\rm Ad}}

\begin{document}
%
%
\title{A splitting result for compact symplectic manifolds}
\author{Lucio Bedulli and Anna Gori}

\address{Lucio Bedulli\\{\it Dipartimento di Matematica U.Dini\\
Viale Morgagni 67 A}, 50134 Firenze\\{\it Italy.}$\,\,\,${\it E-mail address: }
{\tt bedulli@math.unifi.it}}
\address{Anna Gori\\{\it Dipartimento di Matematica e Appl.\ per l'Architettura\\Piazza Ghiberti 27}, 50100 Firenze\\{\it Italy}$\,\,\,${\it E-mail address: }{\tt gori@math.unifi.it}}
\begin{abstract}
We consider compact symplectic manifolds acted on effectively 
by a compact connected Lie group $K$ in a Hamiltonian fashion.
We prove that the squared moment map $||\mu||^2$ is constant  if and
only if $K$ is semisimple and the manifold is $K$-equivariantly
symplectomorphic to a product of a
flag manifold and a compact symplectic manifold which is acted on
trivially by $K$. In the almost-K\"ahler setting the symplectomorphism
turns out to be an isometry.
\end{abstract}
\maketitle
\section{Introduction}
In \cite{GP1} it has been proved that on a compact K\"ahler
manifold acted on effectively by a compact connected Lie group $K$
of isometries in a Hamiltonian fashion, {\em  the squared moment map
$||\mu||^2$ is constant if and only if $K$ is semisimple and the
manifold is biholomorphically and $K$-equivariantly isometric to a
product of a flag manifold and a compact K\"ahler manifold which
is acted on trivially by $K$.}\\ In the present paper we consider
its symplectic analogue; we deal with a $2n$-dimensional symplectic manifold $M$
acted on by a compact connected Lie group $K$, we suppose that the
$K$-action on $M$ is Hamiltonian, i.e. there exists a moment map
$\mu: M \to \k^*$, where $\k$ is the Lie algebra of $K$.
Throughout the following we will denote by $\omega$ the
symplectic form on $M$; moreover Lie groups
and their Lie algebras will be indicated by capital and gothic letters respectively. \\
If we fix an $\Ad(K)$-invariant scalar product $q:=
\langle,\rangle$ on $\k$ and we identify $\k^*$ with $\k$ by means
of $q$, we can think of $\mu$ as a $\k$-valued map; the function
$f\in C^\infty(M)$ defined as $f:= ||\mu||^2$ has been extensively
used in \cite{Kir1} to obtain strong information on the topology of the manifold. \\
Our result is the following
\begin{theorem}\label{splitting} Suppose
 $M$ is a compact symplectic  $K$-Hamiltonian manifold, where $K$ is a
compact  connected  Lie group  acting  effectively  on  $M$. If  $\mu$
is the corresponding moment map, then its squared norm $f=\|\mu\|^2$ is
constant if  and only  if $K$  is semisimple and  the manifold  $M$ is
$K$-equivariantly symplectomorphic  to the product of  a flag manifold
and a compact manifold which is acted on trivially by $K$. Moreover if
on $(M,\omega)$ is given a $K$-invariant $\omega$-compatible almost
complex structure, the symplectomorphism turns out to be an isometry
with respect to the induced Riemannian metric.
\end{theorem}
In the first subsection of the section 2 we give the proof of the
symplectic part of the statement of Theorem \ref{splitting}, while in
the second subsection we deal with the almost-K\"ahler setting. 
\section{Proof of the main result}
\subsection{The symplectic setting}
We will follow the notations as in \cite{Kir1}.
Assume $f$ to be constant, i.e. $\mu$ maps the manifold $M$ into a
sphere. We fix a point $x_o\in M$; we can suppose that $\beta=\mu(x_o)$ lies in the
closure of a Weyl chamber $\Lt_+$, where $\Lt$ denotes the Lie
algebra of a fixed maximal torus in $K$. Since $\mu(M)$ is
connected and $\mu(M)\cap \mathfrak{t}^*_+$ is convex \cite{Kir2},
$\mu(M)$ is a single coadjoint orbit $\mathcal{O}$ $=K/K_{\beta}$,
where $\beta = \mu(x_o)$. \\ We first want to prove that
$\mu^{-1}(\beta)$ is a symplectic submanifold of $(M,\omega)$. Let $\mu_\beta :=
\langle\mu,\beta\rangle$ be the height function relative to
$\beta$ and let $Z_\beta$ denote the intersection of the critical
point set of $\mu_\beta$, with the pre-image of $||\beta||^2$ via
$\mu_\beta$. Since $f$ is constant,
each point in $M$ is critical for $f$, therefore its critical set,
which is proved (Lemma 3.15 \cite{Kir1}) to be
$C_\beta=K\cdot(Z_\beta\cap \mu^{-1}(\beta))$,
 coincides with $M$. In Kirwan's language we can say that the
 stratification $\{S_\beta\}$, induced by $f$, is made of only one stratum. Now,
following the proof of Proposition $2$ in \cite{GP1}, we have that
$Z_\beta = \mu^{-1}(\beta)$. Indeed, if $p\in Z_\beta$, then
$\mu_\beta(p) = ||\beta||^2$ and
\[
||\beta||^2 = \langle\mu(p),\beta\rangle \leq ||\mu(p)||\cdot
||\beta|| \leq ||\beta||^2,
\]
and therefore $\mu(p) = \beta$, i.e.
$p\in \mu^{-1}(\beta)$. Vice versa, if $p\in \mu^{-1}(\beta)$,
then $||\mu(p)||^2$ is a critical value of $f:= ||\mu||^2$ and
therefore $\hat{\beta}_p = 0$, where $\hat\beta$ denotes the
fundamental field on $M$ induced by the element $\beta\in \k$;
moreover $\mu_\beta(p) = ||\beta||^2$, hence $p\in
Z_\beta$. This implies that $\mu^{-1}(\beta)$ is a symplectic
submanifold \cite{GS2} and that $M=C_\beta = K\cdot
\mu^{-1}(\beta).$ As a consequence, at every $y\in \mu^{-1}(\beta)$, the tangent
space to $M$ splits as
\[
T_y M = T_y(K\cdot y) \oplus T_y(\mu^{-1}(\beta)).
\]
Using the fundamental property of the moment map, we can see that the
previous splitting is  symplectic. Indeed if $v\in
T_y(\mu^{-1}(\beta))$ and $X\in \k$ we have
\begin{equation}
\label{Tsplit}
0 = \langle d\mu_y(w),X\rangle = \omega_y(w,{\hat X}_y).
\end{equation}
We have used the fact that $T_y
\mu^{-1}(\beta)=\Ker \,d\mu_y,$ because $\mu^{-1}(\beta)$ is a
submanifold of $M.$ 
\\
Note that every level set of $\mu$  is a symplectic submanifold of $M$,
then we can argue as above to obtain the symplectic splitting
\eqref{Tsplit} at every point of $M$.
This is a consequence of the fact that $\mu$ is $K$-equivariant and it
is constant when passing to the quotient.
\\
Here we show that $K$ is semisimple, i.e. the connected component
$Z$ of its center is trivial.  Fix a $K$-principal point $p\in M$
and observe that the restricted map $\mu:K\cdot p\to \mathcal{O}$ is a
covering because $K\cdot p$ is symplectic; indeed $\Ker \, d\mu,$
restricted to the orbit $K\cdot p,$ is trivial (it is the set
${(T_p\,K\cdot p)}^\omega \cap T_p\,K\cdot p$) hence $\mu$ is a local
diffeomorphism and therefore a covering map since the orbit $K\cdot p$
is compact. 
Since the coadjoint orbits of a compact connected Lie group are
 simply connected, $K\cdot p$ is simply connected too,
 hence $Z$ acts trivially
on it (see e.g.~\cite{B},~p. 224). Since $p$ is principal, $Z$
acts trivially on $M$,
hence $Z$ is trivial, because the $K$-action is effective.\\

We claim that for all $x\in M,$ $K_x=K_{\beta}$. By the
$K$-equivariance of $\mu,$ we have that
$K_x\subseteq K_{\mu(x)}=K_\beta$, which is connected, being the
centralizer of a torus in a compact Lie group
. 
The other inclusion follows from the fact that the
map $\mu:K\cdot p\to \mathcal{O}$ is a covering and both $K\cdot
p$ and $\mathcal{O}$ are simply connected, therefore $\mathfrak
{k}_o=\mathfrak{k}_{\beta}$ and we get the equality.\\
Note that each
$K$-orbit intersects $\mu^{-1}(\beta)$ in a single point:
if there are two points $x$ and $z=k\cdot x$ for $k\in K$ which
lie in $\mu^{-1}(\beta)$, then, by the $K$-equivariance of $\mu,$
we have $k\in K_{\beta}=K_x$; hence $k\cdot x = x$.\\
From this it follows that the map
\[
\varphi: K/K_\beta\times \mu^{-1}(\beta)\to M,\qquad
\varphi(gK_\beta,x) = g\cdot x,
\]
where we identify $K/K_\beta$ with the orbit $K\cdot x_o$, is a
well defined $K$-equivariant diffeomorphism.\\
We also observe that $\mu^{-1}(\mu(x))$ is connected for all $x\in
M$. Moreover all the $K$-orbits are principal since their
stabilizers are all equal to $K_{\beta}$, hence we have that $K_x$
acts trivially on $T_x \mu^{-1}(\mu(x))=(T_x K\cdot
x)^\omega$ for all $x\in M$. \\
We now denote by $\mathcal F$ the foliation given by the
$K$-orbits and by ${\mathcal F}^\omega$ the $\omega$-orthogonal foliation,
so that ${\mathcal F}^\omega_y = T_y(\mu^{-1}(\mu(y)))$. By the same
symbol we denote the corresponding integrable distributions.\\
We claim that $\varphi$ is a symplectomorphism. It is sufficient
to prove it locally, the global property following from the fact that $\varphi$
is a diffeomorphism. Choose coordinates
$x_1,x_2,\ldots,x_{2n}$ in a neighborhood $U$ of a point $p\in
M,$ in such a way that $\{X_1=\frac{\partial}{\partial
  {x_i}},\ldots,X_{2r}=\frac{\partial}{\partial {x_{2r}}}\}$ is a local
frame for $\mathcal {F}$, while $\{X_{2r+1}=\frac{\partial}{\partial
  {x_{2r+1}}},\ldots,X_{2n}=\frac{\partial}{\partial {x_{2n}}}\}$ is a local
frame for $\mathcal {F}^\omega$.
We show that $M$ is locally symplectomorphic to the product
$K/K_\beta\times \mu^{-1}(\beta)$ proving that the following three
conditions are satisfied
\begin{itemize}
\item[(i)] $\omega(X_i,X_j)=0$ for all $i=1,\ldots,2r$ and $j=2r+1,\ldots,2n$.
\item[(ii)] $X_l\omega(X_i,X_j)=0$ for all $l=1\ldots 2r$
and $i,j=2r+1\ldots 2n$
\item[(iii)] $X_m\omega(X_i,X_j)=0$ for all $m=2r+1\ldots 2n$
and $i,j=1\ldots 2r.$
\end{itemize}
We have already proved that $ \omega(v,w)=0$ for all
$v\in\Gamma({\mathcal F})$ and $w \in \Gamma({\mathcal F}^\omega)$
using the fundamental property of $\mu$. Nevertheless this fact can
also be seen as a consequence of the following remark. Let $x \in M$
be arbitrary. Fix $v\in
\Gamma(\mathcal{F}_x^\omega),$ and consider the linear form
$\mathcal{L}_v\in (T_x K\cdot x)^*$ assigning to each 
$w\in\Gamma(\mathcal{F}_x)$ the real number $\omega(v,w).$ We claim
that it is always zero, hence, in particular, $(i)$ follows. 
Since the $K$-action is symplectic and $x$ is a principal point,   
so that the isotropy representation of $K_x$ leaves ${\mathcal
F}^\omega_x$ pointwisely fixed, ${\mathcal L}_v$ is $K_x$-invariant. Indeed:
\[
\mathcal{L}_v({L_k}_* w)=\omega(v,{L_k}_* 
w)=\omega({{L_k}}^{-1}_*(v),w)=\omega(v,w)=\mathcal{L}_v(w)
\]
On the other hand the isotropy representation of $K_x$ on $T_x\,K\cdot
x$ (hence on $(T_x \, K \cdot x)^*$) has no nonzero fixed vector, since
$K_x$ is the
centralizer of some torus in $K$, therefore $\mathcal{L}_v \equiv 0$. To
prove $(ii)$, recall that the closedness of $\omega$ means  
\[
X\omega(Y,Z)+Y\omega(Z,X)+Z\omega(X,Y)+
\]
\[
-\omega([X,Y],Z)+\omega([X,Z],Y)-\omega([Y,Z],X)=0.
\]
for every $X,Y,Z \in \Gamma(TM)$.
Using the coordinates introduced above and $(i)$, the claim follows.
The same argument applies to $(iii).$
\subsection{The almost-K\"ahler setting}
 Suppose now that on $(M,\omega)$ is given a $K$-invariant
 $\omega$-compatible almost complex structure $J.$ Denote by $g$ the
 induced Riemannian metric defined as $g(X,Y):=\omega(X,JY)$ for all
 $X,Y\in \Gamma(TM)$.\\
In the previous part we  proved that $\mu^{-1}(\beta)=Z_\beta.$ In
 \cite{Kir1} (Lemma 4.12) it is shown that $Z_\beta$ is $J$-invariant
 w.r.t. a fixed $\omega$-compatible almost-complex structure. This
 allows us to deduce that the symplectic splitting \eqref{Tsplit} is
 also $g$-orthogonal, hence every $K$-orbit is $J$-invariant.\\
In order to prove that $\varphi$ is a local isometry now it is sufficient
  to show that
 the metric analogues of $(ii)$ and $(iii)$ hold. \\
We use the same notations introduced in the previous subsection. The
 proofs of these conditions are completely analogous, we prove for
 example the first one. 
Suppose  $X_i$ and the pair $X_j,X_l$  are local sections around $p$ of $\mathcal{F}$
 and $\mathcal{F}^\omega$
 respectively. We have that $X_i \, g(X_j,X_l)=X_i \,\omega(X_j,J X_l);$
 since $\mathcal{F}^\omega$ is $J$-invariant, we can express $J X_l$ as
  $\sum_{m=2r+1}^{2n}a_m(x)X_m$, where $a_m$ are $C^\infty$-functions
 on $U$ which do not depend on the first $2r$ variables. Therefore
 the  metric condition $(ii)$  becomes:
\begin{eqnarray*}
X_i\,\omega(X_j,J X_l) & = &
\sum_{m=2r+1}^{2n}X_i (a_m(x))\cdot
 \omega(X_j,X_m)+\\
 & & \sum_{m=2r+1}^{2n}a_m(x)X_i\omega(X_j,X_m),
\end{eqnarray*}
which vanishes because of the symplectic $(ii).$ 
\begin{remark}The compactness assumption on $M$ can be weakened assuming that the moment map
$\mu$ is proper. In this case, as  Knop shows in \cite{K}, we have that $\mu(M)\cap \mathfrak{t_+}^*$ is convex. Therefore $M$ is still
 mapped to a coadjoint orbit, which is compact. Therefore $M$ must be compact.
\end{remark}    

\end{document}